\documentclass[10pt]{article}
\usepackage[dvips]{color}

\newlength{\minitwocolumn}
\font\teneufm=eufm10
\font\seveneufm=eufm7
\font\fiveeufm=eufm5
\newfam\eufmfam
\textfont\eufmfam=\teneufm
\scriptfont\eufmfam=\seveneufm
\scriptscriptfont\eufmfam=\fiveeufm

\makeatletter
\@addtoreset{equation}{section}
\makeatother

\newtheorem{thm}{Theorem}[section]
\newtheorem{prop}[thm]{Proposition}

\title{\Large{\bf Wakimoto realization of \\
the quantum affine superalgebra $U_q(\widehat{sl}(M|N))$}}
\begin{document}
\maketitle

\begin{center}
{TAKEO KOJIMA}
\end{center}

\begin{abstract}
A bosonization of the quantum affine superalgebra $U_q(\widehat{sl}(M|N))$ is presented for an arbitrary level $k \in {\bf C}$.
The Wakimoto realization is given by using $\xi-\eta$ system.
The screening operators that commute with $U_q(\widehat{sl}(M|N))$ are presented for the level $k \neq -M+N$.
New bosonization of the affine superalgebra $\widehat{sl}(M|N)$ is obtained in the limit $q \to 1$.
\end{abstract}

\section{Introduction}
\label{sec:1}

Bosonization is a powerful method to study representation theory 
and its application to mathematical physics \cite{Frenkel}.
Wakimoto realization is the bosonization
that provides a bridge between representation theory of affine algebras and the geometry of the semi-infinite flag manifold.
The Wakimoto realizations have been constructed for the affine Lie algebra $g=(ADE)^{(r)}~(r=1,2)$, $(BCFG)^{(1)}$ and 
$\widehat{sl}(M|N)$, ${osp}(2|2)^{(2)}$, $D(2,1,a)^{(1)}$ \cite{Wakimoto, Feigin-Frenkel1, Feigin-Frenkel2, Ito-Komata, Boer-Feher, Szczesny, Feher-Pusztai, Ding-Gould-Zhang, Yang-Zhang-Liu, Iohara-Koga, Shafiekhani-Chung}.
They have been used to construct correlation functions of WZW models, in the study of Drinfeld-Sokolov reduction and $W$-algebras.
It's nontrivial to give quantum deformation of Wakimoto realization as the same as quantum
Drinfeld-Sokolov reduction and quantum $W$-algebras.
The quantum Wakimoto realizations have been constructed only for $U_q(\widehat{sl}(N))$ and $U_q(\widehat{sl}(2|1))$
\cite{Matsuo, Shiraishi, Awata-Odake-Shiraishi1, Awata-Odake-Shiraishi2, Zhang-Gould}.  
In this paper we study a higher-rank generalization of the previous works
for the quantum affine superalgebra $U_q(\widehat{sl}(2|1))$. 
We give a bosonization of the quantum affine superalgebra $U_q(\widehat{sl}(M|N))$ for an arbitrary level $k \in {\bf C}$,
and give the Wakimoto realization using $\xi-\eta$ system.
We give the screening operators that commute with $U_q(\widehat{sl}(M|N))$ for the level $k \neq -M+N$.
Taking the limit $q \to 1$, we obtain new bosonization of the affine superalgebra $\widehat{sl}(M|N)$.
This paper is a shorter review of the papers \cite{Kojima1, Kojima2, Kojima3, Kojima4}.

\section{Quantum affine superalgebra $U_q(\widehat{sl}(M|N))$}
\label{sec:2}

In this Section we recall the definition of the quantum affine superalgebra $U_q(\widehat{sl}(M|N))$ for $M,N=1,2,3,\cdots$.
Throughout this paper, $q \in {\bf C}$ is assumed to be $0<|q|<1$.
For any integer $n$, define $[n]_q=\frac{q^n-q^{-n}}{q-q^{-1}}$.
We set $\nu_i=+1~(1\leq i \leq M)$, $\nu_i=-1~(M+1\leq i \leq M+N)$ and $\nu_0=-1$.
The Cartan matrix $(A_{i,j})_{0 \leq i, j \leq M+N-1}$ of the affine Lie superalgebra $\widehat{sl}(M|N)$ is
given by 
$$A_{i,j}=(\nu_i+\nu_{i+1})\delta_{i,j}-\nu_i \delta_{i,j+1}-\nu_{i+1}\delta_{i+1,j}.$$

The quantum affine superalgebra $U_q(\widehat{sl}(M|N))$ \cite{Yamane} is
the associative algebra over ${\bf C}$ with the generators
$X_{m}^{\pm, i}~(i=1,2,\cdots,M+N-1, m \in {\bf Z})$,
$H_{n}^i~(i=1,2,\cdots,M+N-1, n \in {\bf Z}_{\neq 0})$,
$H^i~(i=1,2,\cdots,M+N-1)$, and $c$.
The ${\bf Z}_2$-grading of the generators is given by
$p(X_m^{\pm, M})\equiv 1 \pmod{2}$ for $m \in {\bf Z}$ and zero otherwise.
The defining relations of the generators are given as follows.
\begin{eqnarray}
&&c : {\rm central~element},\nonumber\\
&&[H^i,H_m^j]=0,~~~[H_{m}^i,H_{n}^j]=\frac{[A_{i,j}m]_q[cm]_q}{m}\delta_{m+n,0},
\nonumber\\
&&[H^i,X^{\pm,j}(z)]=\pm A_{i,j}X^{\pm,j}(z),
\nonumber\\
&&[H_{m}^i, X^{\pm,j}(z)]=\pm \frac{[A_{i,j}m]_q}{m}q^{\mp \frac{c}{2}|m|} z^m X^{\pm,j}(z),
\nonumber
\\
&&(z_1-q^{\pm A_{i,j}}z_2)
X^{\pm,i}(z_1)X^{\pm,j}(z_2)
=
(q^{\pm A_{j,i}}z_1-z_2)
X^{\pm,j}(z_2)X^{\pm,i}(z_1),~~{\rm for}~|A_{i,j}|\neq 0,
\nonumber
\\
&&
[X^{\pm,i}(z_1), X^{\pm,j}(z_2)]=0~~~{\rm for}~|A_{i,j}|=0,
\nonumber
\\
&&[X^{+,i}(z_1), X^{-,j}(z_2)]
=\frac{\delta_{i,j}}{(q-q^{-1})z_1z_2}
\left(
\delta(q^{c}z_2/z_1)\Psi_+^i(q^{\frac{c}{2}}z_2)-
\delta(q^{-c}z_2/z_1)\Psi_-^i(q^{-\frac{c}{2}}z_2)\right), \nonumber
\\
&& 
[X^{\pm,i}(z_{1}),
[X^{\pm,i}(z_{2}), X^{\pm,j}(z)]_{q^{-1}}]_q+\left(z_1 \leftrightarrow z_2\right)=0
~~~{\rm for}~|A_{i,j}|=1,~i\neq M,
\nonumber\\
&&
[X^{\pm,M}(z_1), [X^{\pm,M+1}(w_1), [X^{\pm,M}(z_2), X^{\pm, M-1}(w_2)]_{q^{-1}}]_q ]
+(z_1 \leftrightarrow z_2)=0,
\nonumber
\end{eqnarray}
where we use
\begin{eqnarray}
[X,Y]_a=XY-(-1)^{p(X)p(Y)}a YX,\nonumber
\end{eqnarray}
for homogeneous elements $X,Y \in U_q(\widehat{sl}(M|N))$.
For simplicity we write $[X,Y]=[X,Y]_1$.
Here we set
$\delta(z)=\sum_{m \in {\bf Z}}z^m$
and the generating functions
\begin{eqnarray}
&&X^{\pm,j}(z)=
\sum_{m \in {\bf Z}}X_{m}^{\pm,j} z^{-m-1},\nonumber\\
&&\Psi_\pm^i(q^{\pm \frac{c}{2}}z)=q^{\pm h_i}
\exp\left(
\pm (q-q^{-1})\sum_{m>0}H_{\pm m}^i z^{\mp m}
\right).\nonumber
\end{eqnarray}

The multiplication rule for the tensor product is ${\bf Z}_2$-graded and is defined for homogeneous elements 
$X_1, X_2, Y_1, Y_2 \in U_q(\widehat{sl}(M|N))$ by
$(X_1 \otimes Y_1) (X_2 \otimes Y_2)=(-1)^{p(Y_1)p(X_2)} (X_1 X_2 \otimes Y_1 Y_2)$,
which extends to inhomogeneous elements through linearity.

Let $\bar{\alpha}_i$, $\bar{\Lambda}_i$ $(1\leq i \leq M+N-1)$ be the classical simple roots,
the classical fundamental weights, respectively.
Let $(\cdot|\cdot)$ be the symmetric bilinear form satisfying
$(\bar{\alpha}_i|\bar{\alpha}_j)=A_{i,j}$ and $(\bar{\Lambda}_i|\bar{\alpha}_j)=\delta_{i,j}$ for $1\leq i,j \leq M+N-1$.
Let us introduce the affine weight $\Lambda_0$ and the null root $\delta$ 
satisfying $(\Lambda_0|\Lambda_0)=(\delta|\delta)=0$, $(\Lambda_0|\delta)=1$,
and $(\Lambda_0|\bar{\alpha}_i)=(\Lambda_0|\bar{\Lambda}_i)=0$ for $1\leq i \leq M+N-1$.
The other affine weights and the affine roots are given by
$\Lambda_i=\bar{\Lambda}_i+\Lambda_0$, $\alpha_i=\bar{\alpha}_i$ for $1 \leq i \leq M+N-1$,
and $\alpha_0=\delta-\sum_{i=1}^{M+N-1}\alpha_i$.
Let $V(\lambda)$ be the highest-weight module over $U_q(\widehat{sl}(M|N))$ generated by the highest weight vector $|\lambda \rangle \neq 0$ such that
\begin{eqnarray}
&&H_m^i |\lambda \rangle=X_m^{\pm,i}|\lambda \rangle=0~~~(m>0),\nonumber
\\
&&X_0^{+,i}|\lambda \rangle=0,~~~H^i|\lambda \rangle=l_i |\lambda \rangle,\nonumber
\end{eqnarray}
where the classical part of the highest weight is $\bar{\lambda}=\sum_{i=1}^{M+N-1} l_i \bar{\Lambda}_i$.

\section{Bosonization of $U_q(\widehat{sl}(M|N))$}

In this Section we give a bosonization of $U_q(\widehat{sl}(M|N))$ for an arbitrary level $k \in {\bf C}$.

\subsection{Boson}

We introduce bosons 
$a_m^i~(m \in {\bf Z}, 1\leq i \leq M+N-1)$,
$b_m^{i,j}~(m \in {\bf Z}, 1\leq i<j \leq M+N)$,
$c_m^{i,j}~(m \in {\bf Z}, 1\leq i<j \leq M+N)$,
and zero mode operators
$Q_a^i~(1\leq i \leq M+N-1)$,
$Q_b^{i,j}~(1\leq i<j \leq M+N)$,
$Q_c^{i,j}~(1\leq i<j \leq M+N)$.
Their commutation relations are
\begin{eqnarray}
&&[a_m^i,a_n^j]=\frac{1}{m}[(k+g)m]_q[A_{i,j}m]_q\delta_{m+n,0},
~~~[a_0^i,Q_a^j]=(k+g)A_{i,j},\nonumber
\\
&&[b_m^{i,j},b_n^{i',j'}]=-\nu_i \nu_j \frac{1}{m} [m]_q^2 \delta_{i,i'}\delta_{j,j'}\delta_{m+n,0},
~~~[b_0^{i,j},Q_b^{i',j'}]=-\nu_i \nu_j \delta_{i,i'} \delta_{j,j'},
\nonumber
\\
&&[c_m^{i,j},c_n^{i',j'}]=\nu_i \nu_j \frac{1}{m} [m]_q^2 \delta_{i,i'}\delta_{j,j'}\delta_{m+n,0},
~~~[c_0^{i,j},Q_c^{i',j'}]=\nu_i \nu_j \delta_{i,i'} \delta_{j,j'},
\nonumber\\
&&[Q_b^{i,j},Q_b^{i',j'}]=\pi \sqrt{-1}~~~(\nu_i \nu_j=\nu_{i'}\nu_{j'}=-1).\nonumber
\end{eqnarray}
The remaining commutators vanish.
Here $g=M-N$
stands for the dual Coxeter number.
We define free boson fields $b_\pm^{i,j}(z), b^{i,j}(z)$ as follows.
\begin{eqnarray}
&&b_\pm^{i,j}(z)=
\pm (q-q^{-1}) \sum_{m>0} b_{\pm m}^{i,j} z^{\mp m}\pm b_0^{i,j} {\rm log}q,\nonumber
\\
&&b^{i,j}(z)=-\sum_{m \neq 0}\frac{b_m^{i,j}}{[m]_q}z^{-m}+Q_b^{i,j}+b_0^{i,j}{\rm log}z.\nonumber
\end{eqnarray}
Free boson fields $a_\pm^i(z), c^{i,j}(z)$ are defined in the same way.
We define free boson fields $(\Delta^{\varepsilon}_{L} b_\pm^{i,j})(z),
(\Delta^{\varepsilon}_{R} b_\pm^{i,j})(z)$ $(\varepsilon=\pm,0)$ as follows.
\begin{eqnarray}
&&(\Delta^{\varepsilon}_L b_\pm^{i.j})(z)=
\left\{\begin{array}{cc}
b_\pm^{i+1,j}(q^{\varepsilon}z)-b_\pm^{i,j}(z)& (\varepsilon=\pm),
\nonumber\\
b_\pm^{i+1,j}(z)+b_\pm^{i,j}(z)& (\varepsilon=0),
\end{array}\right.
\\
&&(\Delta^{\varepsilon}_R b_\pm^{i.j})(z)=\left\{\begin{array}{cc}
b_\pm^{i,j+1}(q^{\varepsilon}z)-b_\pm^{i,j}(z)& (\varepsilon=\pm),\\
b_\pm^{i,j+1}(z)+b_\pm^{i,j}(z)& (\varepsilon=0).
\end{array}\right.
\nonumber
\end{eqnarray}
We define
free boson fields with parameters $L_1,\cdots,L_r, M_1,\cdots, M_r, \alpha$ as follows.
\begin{eqnarray}
&&
\left(\frac{L_1}{M_1}\frac{L_2}{M_2}\cdots \frac{L_r}{M_r}~a^i\right)(z;\alpha)\nonumber
\\
&=&-\sum_{m \neq 0}
\frac{[L_1 m]_q [L_2 m]_q \cdots [L_r m]_q}{[M_1 m]_q [M_2m]_q \cdots [M_r m]_q}\frac{a_m^i}{[m]_q}q^{-\alpha|m|}z^{-m}
+\frac{L_1 L_2 \cdots L_r}{M_1 M_2 \cdots M_r}(Q_a^i+a_0^i {\rm log}z).\nonumber
\end{eqnarray} 
Normal ordering rules are defined as follows.
\begin{eqnarray}
&&:b_m^{i,j} b_n^{i',j'}:=:b_n^{i',j'} b_m^{i,j}:=\left\{
\begin{array}{cc}
b_m^{i,j} b_n^{i',j'}& (m<0),\\
b_n^{i',j'} b_m^{i,j}& (m>0),
\end{array}
\right.\nonumber
\\
&&:Q_b^{i,j} Q_b^{i',j'}:=:Q_b^{i',j'} Q_b^{i,j}:=Q_b^{i,j} Q_b^{i',j'}~~~(i>i'~~{\rm or}~~i=i',j>j').
\nonumber
\end{eqnarray}
Normal ordering rules of $a_m^i$, $c_m^{i,j}$ and $Q_c^{i,j}$ are defined in the same way.

\subsection{Bosonization}

We define bosonic operators $\Psi_\pm^i(z)~(1\leq i \leq M+N-1)$ as follows.
\begin{eqnarray}
\Psi_\pm^i(q^{\pm \frac{k}{2}}z)
&=&
:e^{a_\pm^i(q^{\pm \frac{g}{2}}z)+\sum_{l=1}^i (\Delta_R^{\mp}b_\pm^{l,i})(q^{\pm(\frac{k}{2}+l)}z)
-\sum_{l=i+1}^M (\Delta_L^{\mp}b_\pm^{i,l})(q^{\pm(\frac{k}{2}+l)}z)}\nonumber\\
&\times& 
e^{-\sum_{l=M+1}^{M+N}(\Delta_L^{\mp} b_\pm^{i,l})(q^{\pm(\frac{k}{2}+2M+1-l)}z)}:~~~(1\leq i \leq M-1),
\label{def:psi1}
\\
\Psi_\pm^M(q^{\pm \frac{k}{2}}z)
&=&
:e^{a_\pm^M(q^{\pm \frac{g}{2}}z)-\sum_{l=1}^{M-1} (\Delta_R^0 b_\pm^{l,M})(q^{\pm(\frac{k}{2}+l)}z)+\sum_{l=M+2}^{M+N} (\Delta_L^0 b_\pm^{M,l})(q^{\pm(\frac{k}{2}+2M+1-l)}z)}:,\nonumber
\\
\label{def:psi2}
\\
\Psi_\pm^i(q^{\pm \frac{k}{2}}z)
&=&:e^{a_\pm^i(q^{\pm \frac{g}{2}}z)
-\sum_{l=1}^M (\Delta_R^\pm b_\pm^{l,i})(q^{\pm (\frac{k}{2}+l-1)}z)-\sum_{l=M+1}^i (\Delta_R^\pm b_\pm^{l,i})(q^{\pm(\frac{k}{2}+2M-l)}z)}\nonumber\\
&\times&
e^{\sum_{l=i+1}^{M+N}
(\Delta_L^\pm  b_\pm^{i,l})(q^{\pm(\frac{k}{2}+2M-l)}z)}:~~~(M+1\leq i \leq M+N-1).
\label{def:psi3}
\end{eqnarray}
We define bosonic operators $X^{\pm,i}(z)~(1\leq i \leq M+N-1)$ as follows.
\begin{eqnarray}
X^{+,i}(z)&=&\sum_{j=1}^i\frac{c_{i,j}}{(q-q^{-1})z}(E_{i,j}^+(z)-E_{i,j}^-(z))~~~(1\leq i \leq M-1),
\label{def:X^+1}\\
X^{+,M}(z)&=&\sum_{j=1}^M c_{M,j} E_{M,j}(z),
\label{def:X^+2}
\\
X^{+,i}(z)&=&
\sum_{j=1}^M c_{i,j} E_{i,j}(z)+\sum_{j=M+1}^i \frac{c_{i,j}}{(q-q^{-1})z}(E_{i,j}^+(z)-E_{i,j}^-(z))\nonumber\\
&&~~~~~(M+1\leq i \leq M+N-1),
\label{def:X^+3}\\
X^{-,i}(z)&=&
\sum_{j=1}^{i-1}
\frac{d_{i,j}^1}{(q-q^{-1})z}(F_{i,j}^{1,-}(z)-F_{i,j}^{1,+}(z))
+\frac{d_{i,i}^2}{(q-q^{-1})z}(F_{i,i}^{2,-}(z)-F_{i,i}^{2,+}(z))\nonumber\\
&+&\sum_{j=i+2}^M \frac{d_{i,j}^3}{(q-q^{-1})z}(F_{i,j}^{3,-}(z)-F_{i,j}^{3,+}(z))+\sum_{j=M+1}^{M+N} d_{i,j}^3 F_{i,j}^3(z)\nonumber\\
&&~~~~~(1\leq i \leq M-1),
\label{def:X^-1}
\\
X^{-,M}(z)
&=&\sum_{j=1}^{M-1}\frac{d_{M,j}^1}{(q-q^{-1})z}(F_{M,j}^{1,-}(z)-F_{M,j}^{1,+}(z))+\frac{d_{M,M}^2}{(q-q^{-1})z}(F_{M,M}^{2,-}(z)-F_{M,M}^{2,+}(z))\nonumber\\
&+&\sum_{j=M+2}^{M+N}\frac{d_{M,j}^3}{(q-q^{-1})z}(F_{M,j}^{3,-}(z)-F_{M,j}^{3,+}(z)),
\label{def:X^-2}
\\
X^{-,i}(z)
&=&\sum_{j=1}^M d_{i,j}^1 F_{i,j}^{1}(z)
+\sum_{j=M+1}^{i-1}\frac{d_{i,j}^1}{(q-q^{-1})z}(F_{i,j}^{1,-}(z)-F_{i,j}^{1,+}(z))\nonumber\\
&&+
\frac{d_{i,i}^2}{(q-q^{-1})z}(F_{i,i}^{2,-}(z)-F_{i,i}^{2,+}(z))
+\sum_{j=i+2}^{M+N}\frac{d_{i,j}^3}{(q-q^{-1})z}(F_{i,j}^{3,-}(z)-F_{i,j}^{3,+}(z))\nonumber\\
&&~~~~~(M+1 \leq i \leq M+N-1).
\label{def:X^-3}
\end{eqnarray}
We set $E_{i,j}^\pm(z)$ as follows.
\begin{eqnarray}
E_{i,j}^\pm (z)&=&:e^{(b+c)^{j,i}(q^{j-1}z)+b_\pm^{j,i+1}(q^{j-1}z)-(b+c)^{j,i+1}(q^{j-1\pm 1}z)
+\sum_{l=1}^{j-1}(\Delta_R^- b_+^{l,i})(q^lz)}:\nonumber\\
&&~~~~~~~(1\leq j<i \leq M-1),\nonumber
\\
E_{i,i}^\pm (z)
&=&:e^{b_\pm^{i,i+1}(q^{i-1}z)-(b+c)^{i,i+1}(q^{i-1\pm 1}z)+\sum_{l=1}^{i-1}
(\Delta_R^- b_+^{l,i})(q^lz)}:~~~(1\leq j<i \leq M-1),
\nonumber\\
E_{i,i}^\pm (z)
&=&:e^{-b_\pm ^{i,i+1}(q^{2M+1-i}z)-(b+c)^{i,i+1}(q^{2M+1 \mp 1-i}z)}\nonumber\\
&\times&
e^{-\sum_{l=1}^M(\Delta_R^+ b_+^{l,i})(q^{l-1}z)
-\sum_{l=M+1}^{i-1}(\Delta_R^+ b_+^{l,i})(q^{2M-l}z)}:~~~(M+1\leq i \leq M+N-1),\nonumber
\\
E_{i,j}^\pm (z)&=&:e^{(b+c)^{j,i}(q^{2M+1-j}z)-b_\pm^{j,i+1}(q^{2M+1-j}z)-(b+c)^{j,i+1}(q^{2M+1 \mp 1-j}z)}\nonumber\\
&\times&
e^{-\sum_{l=1}^M(\Delta_R^+ b_+^{l,i})(q^{l-1}z)-\sum_{l=M+1}^{j-1}(\Delta_R^+ b_+^{l,i})(q^{2M-l}z)}:~~~(M+1\leq j<i\leq M+N-1).\nonumber
\end{eqnarray}
We set $E_{i,j}(z)$ as follows.
\begin{eqnarray}
E_{M,j}(z)
&=&:e^{(b+c)^{j,M}(q^{j-1}z)+b^{j,M+1}(q^{j-1}z)-\sum_{l=1}^{j-1}(\Delta_R^0 b_+^{l,M})(q^lz)}:~~~(1\leq j \leq M-1),\nonumber
\\
E_{M,M}(z)
&=&:e^{b^{M,M+1}(q^{M-1}z)-\sum_{l=1}^{M-1}(\Delta_R^0 b_+^{l,M})(q^lz)}:~~~(1\leq j \leq M-1),
\nonumber\\
E_{i,j}(z)
&=&:e^{b_+^{j,i}(q^{j-1}z)-b^{j,i}(q^jz)+b^{j,i+1}(q^{j-1}z)
-\sum_{l=1}^{j-1}(\Delta_R^+ b_+^{l,i})(q^{l-1}z)}:\nonumber\\
&&~~~~~~(M+1\leq i \leq M+N-1,1\leq j \leq M).\nonumber
\end{eqnarray}
We set $F_{i,j}^{1,\pm}(z), F_{i,j}^1(z)$ as follows.
\begin{eqnarray}
F_{i,j}^{1,\pm}(z)
&=&:e^{a_-^i(q^{-\frac{k+g}{2}}z)+(b+c)^{j,i+1}(q^{-k-j}z)-b_\pm^{j,i}(q^{-k-j}z)-(b+c)^{j,i}(q^{-k-j\mp 1}z)}\nonumber\\
&\times& e^{\sum_{l=j+1}^i (\Delta_R^+ b_-^{l,i})(q^{-k-l}z)-\sum_{l=i+1}^M (\Delta_L^+ b_-^{i,l})(q^{-k-l}z)
-\sum_{l=M+1}^{M+N}(\Delta_L^+ b_-^{i,l})(q^{-k-2M-1+l}z)}:
\nonumber\\
&&~~~~~(1\leq j<i\leq M-1),\nonumber\\
F_{M,j}^{1,\pm}(z)&=&:e^{a_-^M(q^{-\frac{k+g}{2}}z)-b_\pm^{j,M}(q^{-k-j}z)-(b+c)^{j,M}(q^{-k-j \mp 1}z)-b_-^{j,M+1}(q^{-k-j}z)-b^{j,M+1}(q^{-k-j+1}z)}
\nonumber\\
&\times& e^{-\sum_{l=j+1}^{M-1}(\Delta_R^0 b_-^{l,M})(q^{-k-l}z)
+\sum_{l=M+2}^{M+N}(\Delta_L^0 b_-^{M,l})(q^{-k-2M-1+l}z)}:~~~(1\leq j \leq M-1),
\nonumber\\
F_{i,j}^1(z)&=&
:e^{a_-^i(q^{-\frac{k+g}{2}}z)-b_-^{j,i+1}(q^{-k-j}z)-b^{j,i+1}(q^{-k-j+1}z)+b^{j,i}(q^{-k-j}z)
-\sum_{l=j+1}^M(\Delta_R^- b_-^{l,i})(q^{-k-l+1}z)}
\nonumber\\
&\times&
e^{-\sum_{l=M+1}^i (\Delta_R^- b_-^{l,i})(q^{-k-2M+l}z)
+\sum_{l=i+1}^{M+N}(\Delta_L^- b_-^{i,l})(q^{-k-2M+l}z)}:\nonumber\\
&&~~~~~(M+1\leq i \leq M+N-1, 1\leq j \leq M),
\nonumber\\
F_{i,j}^{1,\pm}(z)&=&:e^{a_-^i(q^{-\frac{k+g}{2}}z)+(b+c)^{j,i+1}(q^{-k-2M+j}z)
+b_\pm^{j,i}(q^{-k-2M+j}z)-(b+c)^{j,i}(q^{-k-2M\pm 1+j}z)}\nonumber
\\
&\times& e^{-\sum_{l=j+1}^i(\Delta_R^- b_-^{l,i})(q^{-k-2M+l}z)
+\sum_{l=i+1}^{M+N}(\Delta_L^- b_-^{i,l})(q^{-k-2M+l}z)}:\nonumber\\
&&~~~~~~(M+1\leq j<i \leq M+N-1).\nonumber
\end{eqnarray}
We set $F_{i,i}^{2,\pm}(z)$ as follows.
\begin{eqnarray}
F_{i,i}^{2,\pm}(z)&=&:e^{a_\pm^i(q^{\pm \frac{k+g}{2}}z)
+b_\pm^{i,i+1}(q^{\pm(k+i+1)}z)+(b+c)^{i,i+1}(q^{\pm(k+i)}z)}
\nonumber\\
&\times&
e^{-\sum_{l=i+2}^M(\Delta_L^\mp b_\pm^{i,l})(q^{\pm (k+l)}z)
-\sum_{l=M+1}^{M+N}(\Delta_L^\mp b_\pm^{i,l})(q^{\pm (k+2M+1-l)}z)}:~~~(1\leq i \leq M-1),
\nonumber
\\
F_{M,M}^{2,\pm}(z)
&=&:e^{a_\pm^M(q^{\pm \frac{k+g}{2}}z)-b^{M,M+1}(q^{\pm(k+M-1)}z)
+\sum_{l=M+2}^{M+N}(\Delta_L^0 b_\pm^{M,l})(q^{\pm(k+2M+1-l)}z)}:,
\nonumber\\
F_{i,i}^{2,\pm}(z)&=&:e^{a_\pm^i(q^{\pm \frac{k+g}{2}}z)
-b_\pm^{i,i+1}(q^{\pm(k+2M-1-i)}z)+(b+c)^{i,i+1}(q^{\pm (k+2M-i)}z)}\nonumber\\
&\times& 
e^{\sum_{l=i+2}^{M+N}(\Delta_L^\pm b_\pm^{i,l})(q^{\pm(k+2M-l)}z)}:~~~(M+1\leq i \leq M+N-1).\nonumber
\end{eqnarray}
We set $F_{i,j}^{3,\pm}(z), F_{i,j}^3(z)$ as follows.
\begin{eqnarray}
F_{i,j}^{3,\pm}(z)&=&:e^{a_+^i(q^{\frac{k+g}{2}}z)+
(b+c)^{i,j}(q^{k+j-1}z)
+b_\pm^{i+1,j}(q^{k+j-1}z)-(b+c)^{i+1,j}(q^{k-1\pm 1+j}z)}\nonumber\\
&\times&e^{-\sum_{l=j}^M(\Delta_L^- b_+^{i,l})(q^{k+l}z)
-\sum_{l=M+1}^{M+N}(\Delta_L^- b_+^{i,l})(q^{k+2M+1-l}z)}:~~~(1\leq i<j \leq M-1),
\nonumber\\
F_{i,j}^{3}(z)&=&:e^{a_+^i(q^{\frac{k+g}{2}}z)-b^{i,j}(q^{k+2M-j}z)-b_+^{i+1,j}(q^{k+2M-j}z)+b^{i+1,j}(q^{k+2M+1-j}z)}\nonumber\\
&\times&
e^{-\sum_{l=j+1}^{M+N}(\Delta_L^- b_+^{i,l})(q^{k+2M+1-l}z)}:~~~(1\leq i \leq M-1, M+1\leq j \leq M+N),
\nonumber\\
F_{M,j}^{3,\pm}(z)&=&
:e^{a_+^M(q^{\frac{k+g}{2}}z)-b^{M,j}(q^{k+2M-j}z)-b_\pm^{M+1,j}(q^{k+2M+1-j}z)
-(b+c)^{M+1,j}(q^{k+2M+1 \mp 1-j}z)}\nonumber\\
&\times& e^{b_+^{M+1,j}(q^{k+2M+1-j}z)+\sum_{l=j+1}^{M+N}(\Delta_L^0 b_+^{M,l})(q^{k+2M+1-l}z)}:~~(M+2 \leq j \leq M+N),
\nonumber\\
F_{i,j}^{3,\pm}(z)&=&:e^{a_+^i(q^{\frac{k+g}{2}}z)+(b+c)^{i,j}(q^{k+2M+1-j}z)
-b_\pm^{i+1,j}(q^{k+2M+1-j}z)-(b+c)^{i+1,j}(q^{k+2M+1 \mp1-j}z)}\nonumber\\
&\times& e^{\sum_{l=j+1}^{M+N}(\Delta_L^+ b_+^{i,l})(q^{k+2M-l}z)}:~~~(M+1\leq i<j-1\leq M+N-1).
\nonumber
\end{eqnarray}
The coefficients $c_{i,j} \in {\bf C}$ and 
$d_{i,j}^1, d_{i,i}^2, d_{i,j}^3 \in {\bf C}$ satisfy the following conditions.
\begin{eqnarray}
d_{i,j}^1&=&\nu_{i+1}\frac{1}{c_{i,j}}\times
\left\{
\begin{array}{cc}
1& (1 \leq i \leq M-1, 1 \leq j \leq i-1),
\\
q^{j-1} & (i=M, 1\leq j \leq M-1),
\\
q^{-k-1} & (M+1 \leq i \leq M+N-1, 1 \leq j \leq M),
\\
1 & (M+1 \leq i \leq M+N-1, M+1 \leq j \leq i-1),
\end{array}
\right.\nonumber
\\
d_{i,i}^2&=&
\nu_{i+1}
\frac{1}{c_{i,i}}
\times
\left\{\begin{array}{cc}
1 & (1\leq i \neq M \leq M+N-1),\\
q^{M-1} & (i=M),
\end{array}
\right.\nonumber
\\
d_{i,j}^3
&=&\nu_{i+1}
\frac{1}{c_{i,i}}\prod_{l=1}^{j-i-1}\frac{c_{i+l,i+1}}{c_{i+l,i}}
\times \left\{
\begin{array}{cc}
1& (1 \leq i \leq M-1, i+2 \leq j \leq M),\\
q^{k+3M+1-2j}& (1\leq i \leq M-1, M+1 \leq j \leq M+N),\\
q^{(M-1)(j-M)}& (i=M, M+2 \leq j \leq M+N),\\
1 & (M+1 \leq i \leq M+N-1, i+2 \leq j \leq M+N).
\end{array}\right.
\nonumber
\end{eqnarray}

\begin{thm}~~~The bosonic operators
$\Psi_\pm^i(z)$ defined in (\ref{def:psi1})-(\ref{def:psi3}), and
$X^{\pm,i}(z)$ defined in (\ref{def:X^+1})-(\ref{def:X^+3}) 
and (\ref{def:X^-1})-(\ref{def:X^-3}) satisfy the defining relations of the quantum affine superalgebra $U_q(\widehat{sl}(M|N))$ 
with the central element $c=k \in {\bf C}$.
\label{thm:1}
\end{thm}

\subsection{Wakimoto realization}

In this Section we introduce the $\xi-\eta$ system and give the Wakimoto realization.
We set the boson Fock space ${F}(p_a,p_b,p_c)$ as follows.
The vacuum state $|0\rangle \neq 0$ is defined by
$
a_m^i|0\rangle=
b_m^{i,j}|0\rangle=
c_m^{i,j}|0\rangle=0~(m \geq 0)$.
Let $|p_a,p_b,p_c\rangle$ be
\begin{eqnarray}
&&|p_a,p_b,p_c\rangle
\nonumber\\
&=&\exp\left(
\sum_{i,j=1}^{M+N-1}\frac{(A^{-1})_{i,j}}{k+g}p_a^i Q_a^i-\sum_{1\leq i<j \leq M+N}\nu_i \nu_j p_b^{i,j}Q_b^{i,j}+\sum_{1\leq i<j \leq M+N \atop{\nu_i \nu_j=+1}}p_c^{i,j}Q_c^{i,j}\right)|0\rangle,\nonumber
\end{eqnarray}
then $|p_a, p_b, p_c\rangle$
is the highest weight state of the boson Fock space ${F}(p_a,p_b,p_c)$.
The boson Fock space $F(p_a,p_b,p_c)$ is generated by the bosons $a_m^i, b_m^{i,j}, c_m^{i,j}$ on the highest weight state
$|p_a,p_b,p_c\rangle$.
We set the space $F(p_a)$ by
\begin{eqnarray}
F(p_a)=\bigoplus_{p_b^{i,j}=-p_c^{i,j} \in {\bf Z}~(\nu_i\nu_j=+)
\atop{p_b^{i,j}}\in {\bf Z}~(\nu_i\nu_j=-)}F(p_a,p_b,p_c).\nonumber
\end{eqnarray}
Here we impose the restriction
$p_b^{i,j}=-p_c^{i,j}~(\nu_i \nu_j=+)$, because $X_m^{\pm,i}$ change $Q_b^{i,j}+Q_c^{i,j}$.
$F(p_a)$ is $U_q(\widehat{sl}(M|N))$-module.
Let $|\lambda \rangle=|p_a,0,0\rangle$ where $p_a^i=l_i~(1\leq i \leq M+N-1)$.
The generators $H^i, H_m^i, X_m^{\pm, i}$ act on $|\lambda\rangle$ as follows.
\begin{eqnarray}
&&H_m^i |\lambda \rangle=X_m^{\pm,i}|\lambda \rangle=0~~~(m>0),\nonumber\\
&&X_0^{+,i}|\lambda \rangle=0,~~~H^i|\lambda \rangle=l_i |\lambda \rangle.\nonumber
\end{eqnarray}
We have the level-$k$ highest weight module $V(\lambda)$ of $U_q(\widehat{sl}(M|N))$.
\begin{eqnarray}
V(\lambda) \subset F(p_a).\nonumber
\end{eqnarray}
Here the classical part of the highest weight is $\bar{\lambda}=\sum_{i=1}^{M+N-1} l_i \bar{\Lambda}_i$.

We introduce the $\xi-\eta$ system
We set bosonic operators $\xi_m^{i,j}, \eta_m^{i,j}~(\nu_i \nu_j=+1, 1\leq i<j \leq M+N)$ as follows.
\begin{eqnarray}
\eta^{i,j}(z)=\sum_{m \in {\bf Z}}\eta_m^{i,j} z^{-m-1}=:e^{c^{i,j}(z)}:,~~~
\xi^{i,j}(z)=\sum_{m \in {\bf Z}}\xi_m^{i,j}z^{-m}=:e^{-c^{i,j}(z)}:.\nonumber
\end{eqnarray}
Fourier components 
\begin{eqnarray}
\eta_m^{i,j}=\oint \frac{dz}{2\pi \sqrt{-1}}z^m \eta^{i,j}(z),~~~
\xi_m^{i,j}=\oint \frac{dz}{2\pi \sqrt{-1}}z^{m-1}\xi^{i,j}(z)\nonumber
\end{eqnarray}
are well-defined on the module $F(p_a)$.
The ${\bf Z}_2$-grading is given by
$p(\xi_m^{i,j})=p(\eta_m^{i,j})=+1$.
We have direct sum decomposition.
\begin{eqnarray}
F(p_a)=\eta_0^{i,j}\xi_0^{i,j} F(p_a) \oplus \xi_0^{i,j} \eta_0^{i,j} F(p_a),\nonumber
\end{eqnarray}
where
${\rm Ker}(\eta_0^{i,j})=\eta_0^{i,j} \xi_0^{i,j} F(p_a)$,
${\rm Coker}(\eta_0^{i,j})=\xi_0^{i,j}\eta_0^{i,j} F(p_a)$.
We set
\begin{eqnarray}
\eta_0=\prod_{1\leq i<j \leq M+N\atop{\nu_i\nu_j=+1}}\eta_0^{i,j},~~~
\xi_0=\prod_{1\leq i<j \leq M+N\atop{\nu_i\nu_j=+1}}\xi_0^{i,j}.\nonumber
\end{eqnarray}
We introduce the subspace ${\cal F}(p_a)$ by
\begin{eqnarray}
{\cal F}(p_a)=\eta_0\xi_0 F(p_a).\nonumber
\end{eqnarray}
The operators $\eta_0^{i,j}, \xi_0^{i,j}$ commute with $X^{\pm,i}(z), \Psi_\pm^i(z)$ up to sign $\pm 1$.

\begin{prop}~~~${\cal F}(p_a)$ is $U_q(\widehat{sl}(M|N))$-module.
\end{prop}
We call ${\cal F}(p_a)$ the Wakimoto realization of $U_q(\widehat{sl}(M|N))$.

\section{Screening operator}

In this Section 
we give the screening operators $Q_i~(1\leq i \leq M+N-1)$
that commute with $U_q(\widehat{sl}(M|N))$ for the level $c=k \neq -g$.
We define bosonic operators $S_i(z)~(1\leq i \leq M+N-1)$ that we call the screening currents as follows.
\begin{eqnarray}
S_i(z)&=&\sum_{j=i+1}^M \frac{e_{i,j}}{(q-q^{-1})z}(S_{i,j}^-(z)-S_{i,j}^+(z))+\sum_{j=M+1}^{M+N}e_{i,j}S_{i,j}(z)
\nonumber\\
&&~~~~~~(1\leq i \leq M-1),\label{def:S1}
\\
S_M(z)&=&\sum_{j=M+1}^{M+N}e_{M,j} S_{M,j}(z),
\label{def:S2}
\\
S_i(z)&=&\sum_{j=i+1}^{M+N}\frac{e_{i,j}}{(q-q^{-1})z}(S_{i,j}^-(z)-S_{i,j}^+(z))
\nonumber\\
&&~~~~~(M+1\leq i \leq M+N-1).\label{def:S3}
\end{eqnarray}
We set $S_{i,j}^\pm(z)$ as follows.
\begin{eqnarray}
{S}_{i,j}^\pm(z)&=&:e^{
-(\frac{1}{k+g}a^i)(z;\frac{k+g}{2})
+(b+c)^{i+1,j}(q^{M-N-j}z)-b_\pm^{i,j}(q^{M-N-j}z)
-(b+c)^{i,j}(q^{M-N-j \mp 1}z)}
\nonumber\\
&\times&
e^{\sum_{l=j+1}^{M}(\Delta_L^+ b_-^{i,l})(q^{M-N-l}z)+\sum_{l=M+1}^{M+N}(\Delta_L^+ b_-^{i,l})(q^{-M-N+l-1}z)}:
~~~(1\leq i <j \leq M),\nonumber\\
{S}_{i,j}^\pm(z)&=&:e^{
-(\frac{1}{k+g}a^i)(z;\frac{k+g}{2})
+(b+c)^{i+1,j}(q^{-M-N+j}z)+b_\pm^{i,j}(q^{-M-N+j}z)-(b+c)^{i,j}(q^{-M-N+j\pm 1}z)}\nonumber\\
&\times&
e^{-\sum_{l=j+1}^{M+N}(\Delta_L^- b_-^{i,l})(q^{-M-N+l}z)}:~~~(M+1\leq i<j \leq M+N).\nonumber
\end{eqnarray}
We set $S_{i,j}(z)$ as follows.
\begin{eqnarray}
{S}_{i,j}(z)&=&:e^{
-(\frac{1}{k+g}a^i)(z;\frac{k+g}{2})
+b^{i,j}(q^{-M-N+j}z)+b_+^{i+1,j}(q^{-M-N+j}z)-b^{i+1,j}(q^{-M-N+j+1}z)}\nonumber\\
&\times&e^{\sum_{l=j+1}^{M+N}(\Delta_L^+ b_-^{i,l})(q^{-M-N-1+l}z)}:~~~(1\leq i \leq M-1, M+1\leq j \leq M+N),
\nonumber\\
{S}_{M,j}(z)&=&:e^{
-(\frac{1}{k+g}a^i)(z;\frac{k+g}{2})
+(b+c)^{M+1,j}(q^{-M-N+j}z)+b^{M,j}(q^{-M-N+j}z)}\nonumber\\
&\times& 
e^{-\sum_{l=j+1}^{M+N}(\Delta_L^0 b_-^{M,l})(q^{-M-N-1+l}z)}:~~~(M+1\leq j \leq M+N).\nonumber
\end{eqnarray}
Here we set $e_{i,j}$ as follows.
\begin{eqnarray}
e_{i,i+1}&=&\left\{\begin{array}{cc}
1/d_{i,i}^2& (1\leq i \leq M-1),\\
-q^{-N+1}/d_{M,M}^2& (i=M),\\
-1/d_{i,i}^2& (M+1\leq i \leq M+N-1),
\end{array}\right.
\nonumber
\\
e_{i,j}&=&
\left\{\begin{array}{cc}
1/d_{i,j}^3& (1\leq i \leq M-1, i+2 \leq j \leq M),\\
q^{k+1+M-N}/d_{i,j}^3& (1\leq i \leq M-1, M+1\leq j \leq M+N),\\
-q^{j-M-N}/d_{M,j}^3& (i=M, M+2\leq j \leq M+N),\\
-1/d_{i,j}^3& (M+1\leq i \leq M+N-1, i+2\leq j \leq M+N).
\nonumber
\end{array}
\right.
\end{eqnarray}
The ${\bf Z}_2$-grading of the screening current is given by
$p(S_{M,j}(z))\equiv 1 \pmod{2}$ for $M+1 \leq j \leq M+N$ and zero otherwise.
The Jackson integral with parameters $q \in {\bf C}$ and $s \in {\bf C}^*$ is defined by
\begin{eqnarray}
\int_0^{s \infty}f(w)d_q w=s(1-q)\sum_{n \in {\bf Z}}f(sq^n) q^n.\nonumber
\end{eqnarray}
We define the screening operators $Q_i~(1\leq i \leq M+N-1)$ as follows,
when the Jackson integrals are convergent.
\begin{eqnarray}
Q_i=\int_0^{s \infty}S_i(w)d_{q^{2(k+g)}}w. \label{def:S0}
\end{eqnarray} 

\begin{thm}~~~The screening operators $Q_i$ $(1\leq i \leq M+N-1)$ defined in (\ref{def:S1}), (\ref{def:S2}), (\ref{def:S3}), (\ref{def:S0}) 
commute with the quantum affine superalgebra $U_q(\widehat{sl}(M|N))$.
\begin{eqnarray}
[Q_i,U_q(\widehat{sl}(M|N))]=0.\nonumber
\end{eqnarray}
\end{thm}

\section{Limit $q \to 1$}

Bosonization of the affine superalgebra $\widehat{sl}(M|N)$ for an arbitrary level $k$ have been studied in 
\cite{Ding-Gould-Zhang, Yang-Zhang-Liu, Iohara-Koga}.
We obtain new bosonization of the affine superalgebra $\widehat{sl}(M|N)$ in the limit $q\to 1$.

In what follows we set
\begin{eqnarray}
H^i(z)=\sum_{m \in {\bf Z}} H_m^i z^{-m-1}~~~(1\leq i \leq M+N-1).\nonumber
\end{eqnarray}
We set the parameters $c_{i,j}=1$ in (\ref{def:X^+1})-(\ref{def:X^+3}), (\ref{def:X^-1})-(\ref{def:X^-3}), (\ref{def:S1})-(\ref{def:S3}) for simplicity.
In the limit $q\to1$ we introduce operators $\alpha_i(z)~(1\leq i \leq M+N-1)$,
$\beta_{i,j}(z), \widehat{\beta}_{i,j}(z), \gamma_{i,j}(z)~(1\leq i<j \leq M+N, \nu_i\nu_j=+)$,
and $\psi_{i,j}(z), \psi_{i,j}^\dagger(z)~(1\leq i<j \leq M+N, \nu_i\nu_j=-)$ as follows.
\begin{eqnarray}
&&
\alpha_i(z)=\partial_z \left(a^i(z)\right),
~~~\gamma_{i,j}(z)=:e^{(b+c)^{i,j}(z)}:,
\nonumber\\
&&
\beta_{i,j}(z)=:\partial_z \left(e^{-c^{i,j}(z)}\right) e^{-b^{i,j}(z)}:,~
\widehat{\beta}_{i,j}(z)=:\partial_z \left(e^{-b^{i,j}(z)}\right) e^{-c^{i,j}(z)}:,
\nonumber\\
&&
\psi_{i,j}(z)=:e^{b^{i,j}(z)}:,~~~\psi_{i,j}^\dagger (z)=:e^{-b^{i,j}(z)}:.
\nonumber
\end{eqnarray}
They satisfy the following relations.
\begin{eqnarray}
&&
\alpha_i(z)\alpha_j(w)=\frac{(k+g)A_{i,j}}{(z-w)^2}+\cdots,\nonumber
\\
&&
\beta_{i,j}(z)\gamma_{i',j'}(w)=\frac{\delta_{i,i'}\delta_{j,j'}}{z-w}+\cdots,
~~~
\gamma_{i,j}(z)\beta_{i',j'}(w)=-\frac{\delta_{i,i'}\delta_{j,j'}}{z-w}+\cdots,
\nonumber
\\
&&
\widehat{\beta}_{i,j}(z)\gamma_{i',j'}(w)=-\frac{\delta_{i,i'}\delta_{j,j'}}{z-w}+\cdots,
~~~
\gamma_{i,j}(z)\widehat{\beta}_{i',j'}(w)=\frac{\delta_{i,i'}\delta_{j,j'}}{z-w}+\cdots,
\nonumber
\\
&&
\psi_{i,j}(z)\psi_{i',j'}^\dagger (w)=\frac{\delta_{i,i'}\delta_{j,j'}}{z-w}+\cdots,
~~~
\psi_{i,j}^\dagger(z)\psi_{i',j'}(w)=\frac{\delta_{i,i'}\delta_{j,j'}}{z-w}+\cdots.
\nonumber
\end{eqnarray}
In the limit $q \to 1$ the operators $a_\pm^i(z)$, $b_\pm^{i,j}(z)$, $(\Delta_L^\epsilon b_\pm^{i,j})(z)$ and
$(\Delta_R^\epsilon b_\pm^{i,j})(z)$ disappear.
We obtain the following.
\begin{eqnarray}
H^i(z)&=&\alpha_i(z)+\sum_{j=1}^i:(\widehat{\beta}_{j,i}(z)\gamma_{j,i}(z)-\widehat{\beta}_{j,i+1}(z)\gamma_{j,i+1}(z)):\nonumber\\
&&+\sum_{j=i+1}^M:(\widehat{\beta}_{i+1,j}(z)\gamma_{i+1,j}(z)-\widehat{\beta}_{i,j}(z)\gamma_{i,j}(z)):\nonumber\\
&&+\sum_{j=M+1}^{M+N}:( (\partial_z{\psi}_{i+1,j})(z) \psi_{i+1,j}^\dagger(z)-
(\partial_z {\psi}_{i,j})(z) \psi_{i,j}^\dagger(z)):~~~(1\leq i \leq M-1),\nonumber
\end{eqnarray}
\begin{eqnarray}
H^M(z)&=&\alpha_M(z)+\sum_{j=1}^{M-1}
:((\partial_z {\psi}_{j,M+1})(z)\psi_{j,M+1}^\dagger(z)
+\widehat{\beta}_{j,M}(z)\gamma_{j,M}(z)):\nonumber\\
&&-\sum_{j=M+2}^{M+N}:(\widehat{\beta}_{M+1,j}(z)\gamma_{M+1,j}(z)+(\partial_z {\psi}_{M,j})(z) \psi_{M,j}^\dagger(z)):,
\nonumber
\\
H^i(z)&=&
\alpha_i(z)+\sum_{j=1}^M
:((\partial_z {\psi}_{j,i+1})(z)\psi_{j,i+1}^\dagger(z)-
(\partial_z{\psi}_{j,i})(z)\psi_{j,i}^\dagger(z)):\nonumber\\
&&+\sum_{j=M+1}^i:(\widehat{\beta}_{j,i+1}(z)\gamma_{j,i+1}(z)-\widehat{\beta}_{j,i}(z)\gamma_{j,i}(z)):\nonumber\\
&&+\sum_{j=i+1}^{M+N}
:(\widehat{\beta}_{i,j}(z)\gamma_{i,j}(z)-\widehat{\beta}_{i+1,j}(z)\gamma_{i+1,j}(z)):~~(M+1\leq i \leq M+N-1).
\nonumber
\end{eqnarray}
\begin{eqnarray}
X^{+,i}(z)&=&\sum_{j=1}^i :\beta_{j,i+1}(z)\gamma_{j,i}(z):~~~(1\leq i \leq M-1),
\nonumber\\
X^{+,M}(z)&=&\sum_{j=1}^M :\gamma_{j,M}(z)\psi_{j,M+1}(z):,
\nonumber\\
X^{+,i}(z)&=&\sum_{j=1}^M:\psi_{j,i+1}(z)\psi_{j,i}^\dagger(z):-\sum_{j=M+1}^i
:\beta_{j,i+1}(z)\gamma_{j,i}(z):~(M+1\leq i \leq M+N-1).
\nonumber
\end{eqnarray}
\begin{eqnarray}
X^{-,i}(z)&=&-:\alpha_i(z)\gamma_{i,i+1}(z):-\kappa_i :\partial_z \gamma_{i,i+1}(z):\nonumber\\
&&+\sum_{j=1}^{i-1}:\beta_{j,i}(z)\gamma_{j,i+1}(z):
-\sum_{j=i+2}^M:\beta_{i+1,j}(z)\gamma_{i,j}(z):
-\sum_{j=M+1}^{M+N}:\psi_{i+1,j}(z)\psi_{i,j}^\dagger(z):\nonumber\\
&&+\sum_{j=i+1}^M
:(\widehat{\beta}_{i,j}(z)\gamma_{i,j}(z)
-\widehat{\beta}_{i+1,j}(z)\gamma_{i+1,j}(z))\gamma_{i,i+1}(z):
\nonumber
\\
&&+\sum_{j=M+1}^{M+N}
:((\partial_z {\psi}_{i,j})(z)\psi_{i,j}^\dagger(z)-
(\partial_z{\psi}_{i+1,j})(z)\psi_{i+1,j}^\dagger(z))\gamma_{i,i+1}(z):\nonumber\\
&&~~~~~(1\leq i \leq M-1),
\nonumber\\
X^{-,M}(z)&=&:\alpha_M(z)\psi_{M,M+1}^\dagger(z):+\kappa_M :\partial_z \psi_{M,M+1}^\dagger(z):\nonumber\\
&&-\sum_{j=1}^{M-1}:\beta_{j,M}(z)\psi_{j,M+1}^\dagger(z):
-\sum_{j=M+2}^{M+N}:\beta_{M+1,j}(z)\psi_{M,j}^\dagger(z):\nonumber\\
&&-\sum_{j=M+2}^{M+N}
:(\widehat{\beta}_{M+1,j}(z)\gamma_{M+1,j}^\dagger(z)+
(\partial_z {\psi}_{M,j})(z)\psi_{M,j}^\dagger(z))\psi_{M,M+1}^\dagger(z):,
\nonumber\\
X^{-,i}(z)&=&
:\alpha_i(z)\gamma_{i,i+1}(z):+\kappa_i:\partial_z \gamma_{i,i+1}(z):\nonumber\\
&&-\sum_{j=1}^M:\psi_{j,i}(z)\psi_{j,i+1}^\dagger(z):
+\sum_{j=M+1}^{i-1}:\beta_{j,i}(z)\gamma_{j,i+1}(z):
-\sum_{j=i+2}^{M+N}:\beta_{i+1,j}(z)\gamma_{i,j}(z):\nonumber\\
&&+\sum_{j=i+1}^{M+N}:(\widehat{\beta}_{i,j}(z)\gamma_{i,j}(z)-
\widehat{\beta}_{i+1,j}(z)\gamma_{i+1,j}(z))\gamma_{i,j}(z):\nonumber\\
&&~~~~~(M+1\leq i \leq M+N-1).\nonumber
\end{eqnarray}
Here we have set the coefficients $\kappa_i$ by 
\begin{eqnarray}
\kappa_i=\left\{
\begin{array}{cc}
k+i& (1\leq i \leq M-1)\\
k+M-1& (i=M)\\
k+2M-i& (M+1\leq i \leq M+N-1)
\end{array}\right..\nonumber
\end{eqnarray}
In what follows we assume $k \neq -g$. 
In the limit $q\to 1$ we have the following.
\begin{eqnarray}
S_i(z)&=&\sum_{j=i+1}^M:\tilde{s}_i(z)\beta_{i,j}(z)\gamma_{i+1,j}(z):+
\sum_{j=M+1}^{M+N}:\tilde{s}_i(z)\psi_{i,j}(z)\psi_{i+1,j}^\dagger(z):\nonumber\\
&&~~~~~~(1\leq i \leq M-1),
\nonumber\\
S_M(z)&=&\sum_{j=M+1}^{M+N}:\tilde{s}_M(z)\gamma_{M+1,j}(z)\psi_{M,j}(z):,
\nonumber\\
S_i(z)&=&\sum_{j=i+1}^{M+N}:\tilde{s}_i(z)\beta_{i,j}(z)\gamma_{i+1,j}(z):~~~(M+1\leq i \leq M+N-1).\nonumber
\end{eqnarray}
Here we have set the boson operator $$\tilde{s}_i(z)=:e^{-\left(\frac{1}{k+g}a^i\right)(z;0)}:.$$
Our bosonization is different from \cite{Ding-Gould-Zhang, Yang-Zhang-Liu, Iohara-Koga}.

~\\
{\bf Acknowledgement}~
This work is supported by the Grant-in-Aid for 
Scientific Research {\bf C} (26400105)
from Japan Society for Promotion of Science.
The author would like to thank Professor Michio Jimbo and Professor Vladimir Dobrev for giving advice.
The author would like to thank Professor Zengo Tsuboi, Professor Pascal Baseilhac, Professor Kouichi Takemura and Professor Kenji Iohara for discussion.
The author is thankful for the kind hospitality 
by the organizing committee of 
the 10-th International Symposium "Quantum Theory and Symmetries" (QTS10) and 12-th International Workshop "Lie Theory and Its Applications in Physics" (LT12).

\addcontentsline{toc}{section}{Appendix}

\begin{thebibliography}{99.}%


\bibitem{Frenkel}
E.V.Frenkel, Adv.Math. {\bf 195} (2005) 297-404.

\bibitem{Wakimoto}
M.Wakimoto, Commun.Math.Phys.{\bf 104} (1986) 605-609.

\bibitem{Feigin-Frenkel1}
B.L.Feigin and E.V.Frenkel,
{\it Physics and Mathematics of Strings} (World Scientific, Singapore 1980) 271-316.

\bibitem{Feigin-Frenkel2}
B.L.Feigin and E.V.Frenkel, Commun.Math.Phys. {\bf 128} (1990) 161-189.


\bibitem{Ito-Komata}
K.Ito and S.Komata, Mod.Phys.Lett.{\bf A6} (1991) 581-589.


\bibitem{Boer-Feher}
J.de Boer and L.Feh\'er, Commun.Math.Phys.{\bf 189} (1997) 759-793.

\bibitem{Szczesny}
M.Szczesny, Math.Res.Lett.{\bf 9} (2002) 433-448.


\bibitem{Feher-Pusztai}
L.Feh\'er and B.G.Pusztai, Nucl.Phys.{\bf B674} (2003) 509-532.

\bibitem{Ding-Gould-Zhang}
X.-M.Ding, M.D.Gould and Y.-Z.Zhang, Phys.Lett.{\bf 318} (2003) 354-363.


\bibitem{Yang-Zhang-Liu}
W.-L.Yang,Y.Z.Zhang and X.Liu,
J.Math.Phys.{\bf 48} (2007)  053514 (pp.11).


\bibitem{Iohara-Koga}
K.Iohara and Y.Koga,
Math.Proc.Camb.Phil.Soc.{\bf 132} (2002) 419-433.


\bibitem{Shafiekhani-Chung}
A.Shafiekhani and W.-S.Chung, Mod.Phys.Lett.{\bf A13} (1998) 47-57.


\bibitem{Matsuo}
A.Matsuo, Commun.Math.Phys.{\bf 160} (1994) 33-48.

\bibitem{Shiraishi}
J.Shiraishi, Phys.Lett.{\bf A171} (1992) 243-248.

\bibitem{Awata-Odake-Shiraishi1}
H.Awata, S.Odake and J.Shiraishi, Commun.Math.Phys.{\bf 162} (1994) 61-83.

\bibitem{Awata-Odake-Shiraishi2}
H.Awata, S.Odake and J.Shiraishi,
Lett. Math. Phys.{\bf 42} (1997) 271-279. 

\bibitem{Zhang-Gould}
Y.-Z.Zhang and M.D.Gould, J.Math.Phys.{\bf 41} (2000) 5577-5291.

\bibitem{Kojima1}
T.Kojima, J.Math.Phys.{\bf 53} (2012) 013515 (pp.15).

\bibitem{Kojima2}
T.Kojima, Springer Proceedings {\bf 111} (2013) 263-276.

\bibitem{Kojima3}
T.Kojima, J.Math.Phys.{\bf 53} (2012) 083503 (pp.30).


\bibitem{Kojima4}
T.Kojima, Commun.Math.Phys.{\bf 355} (2017) 603-644.

\bibitem{Yamane}
H.Yamane, Publ.Res.Inst.Math.Sci.{\bf 35} (1999) 321-390.







\end{thebibliography}
\end{document}